\documentclass[preprint,12pt]{article}\usepackage{amsmath,amssymb}\newtheorem{defin}{Definition}[section]
\newtheorem{thm}[defin]{Theorem}\newtheorem{corollary}[defin]{Corollary}\newtheorem{lemma}[defin]{Lemma}
\newtheorem{remark}[defin]{Remark}\newtheorem{example}{Example}
\begin{document}

\begin{center}
{\large\bf Left covariant differential calculi on $\widetilde{\rm GL}_q(2)$}

\vspace*{.3cm}
Salih Celik \\
Department of Mathematics,Yildiz Technical University, DAVUTPASA-Esenler-Istanbul, 34220 TURKEY.
\end{center}

\begin{abstract}
In this work, we introduce the ${\mathbb Z}_3$-graded differential algebra, denoted by $\Omega(\widetilde{\rm GL}_q(2))$, treated as the ${\mathbb Z}_3$-graded quantum de Rham complex of ${\mathbb Z}_3$-graded quantum group $\widetilde{\rm GL}_q(2)$. In this sense, we construct left-covariant differential calculi on the ${\mathbb Z}_3$-graded quantum group $\widetilde{\rm GL}_q(2)$.
\end{abstract}



\section{Introduction}\label{sec1}

Differential geometry of Lie groups plays an important role in the mathematical modeling. When a quantum group (or a Hopf algebra) is considered as a
generalization of the notion of a group, it can be tempting to also generalize the corresponding notions of differential geometry \cite{Brz}.
Differential calculi can be constructed on spaces that are more general than differentiable manifolds. Indeed the algebraic construction of differential calculus in terms of Hopf structures allows to extend the ordinary differential geometric quantities to a variety of interesting spaces that include quantum groups \cite{Woro}. The notions of left-, right- and bi-covariant differential calculus for the matrix quantum groups is introduced by Woronowicz \cite{Woro2}. In the meantime concrete examples of covariant differential calculi on quantum spaces and quantum groups have been constructed by several authors such as \cite{Sch}-\cite{Asc}. In recent years, the ${\mathbb Z}_2$-graded structures \cite{Manin1} has been extensively developed to the ${\mathbb Z}_3$-graded structures (see, for example, \cite{Chung}-\cite{Celik6}), and ${ \mathbb Z}_N$-graded constructions (see, for example, \cite{Dub}, \cite{Kapr}).

The formalism for the construction of the de Rham complex of the quantum groups was proposed by Woronowicz \cite{Woro} and the general theory of this construction was described by Manin \cite{Manin2}. In this paper, differential calculi on the ${\mathbb Z}_3$-graded quantum group $\widetilde{\rm GL}_q(2)$, where $q$ is a primitive cubic root of unity, is presented. Most of this material is new and the Hopf algebraic approach of \cite{Woro} for the construction of differential calculi is adapted to the ${\mathbb Z}_3$-graded quantum group $\widetilde{\rm GL}_q(2)$. Another important property of such a calculus is that the differential {\sf d} satisfies ${\sf d}^3=0$ instead of ${\sf d}^2=0$ such that ${\sf d}^2\ne0$.

\section{${\mathbb Z}_3$-graded concepts}\label{sec2}

In this section, some basic information are presented from ${\mathbb Z}_3$-graded structures in order to ensure that the present article presents an integrity within it and the readers are easily adapted to the subject.

\subsection{${\mathbb Z}_3$-graded algebras}\label{subsec2.1}

The theory of ${\mathbb Z}_3$-graded algebras starts with definition of the concept of ${\mathbb Z}_3$-graded vector space.

\begin{defin} \label{def2.1}
A ${\mathbb Z}_3$-graded vector space $V$ is a vector space over a field ${\mathbb K}$  together with three subspaces $V_0$, $V_1$ and $V_2$ of $V$ such that
$V = V_0 \oplus V_1\oplus V_2$.
\end{defin}
Each subspace $V_i$ is called the {\it i-grade part} of $V$, and its elements are of grade $i$. The grade of an element $v\in V$ is denoted by $p(v)$ and is equal to 0, 1 or 2. All elements of $V$ are collectively said to be  homogeneous.

\begin{defin} \label{def2.2}
Let $L:V\longrightarrow W$ be a linear map of ${\mathbb Z}_3$-graded vector spaces. Then $L$ is called a ${\mathbb Z}_3$-graded vector space homomorphism if it satisfies
$$p(L(v))= p(L) + p(v) \quad ({\rm mod} \, 3) $$
for all $v\in V$.
\end{defin}

\begin{defin} \label{def2.3}
An algebra ${\cal A}$ over ${\mathbb K}$ is called a ${\mathbb Z}_3$-graded algebra if it is a ${\mathbb Z}_3$-graded vector space over ${\mathbb K}$, with a bilinear map ${\cal A}\times{\cal A}\longrightarrow {\cal A}$ such that ${\cal A}_i \cdot{\cal A}_j \subset {\cal A}_{i+j}$ for $i,j=0,1,2$.
\end{defin}

\begin{defin} \label{def2.4}
Let ${\cal A}$ and ${\cal B}$ be ${\mathbb Z}_3$-graded algebras and $L:{\cal A}\longrightarrow {\cal B}$ be a map of definite degree. If it is a
${\mathbb Z}_3$-graded vector space homomorphism and
$$L(ab) = q^{p(a)p(L(b))} L(a) L(b)$$
for all $a, b \in {\cal A}$, where $q$ is a cubic root of unity, then $L$ is called a ${\mathbb Z}_3$-graded algebra homomorphism.
\end{defin}

\begin{defin} \label{def2.5}
Let ${\cal A}$ be a ${\mathbb Z}_3$-graded algebra and the map $L:{\cal A}\longrightarrow {\cal A}$ be a ${\mathbb Z}_3$-graded vector space homomorphism. If it satisfies the ${\mathbb Z}_3$-graded Leibniz rule
$$L(ab) = L(a)b + q^{p(L)p(a)} a L(b), \qquad \forall a,b \in {\cal A}$$
where $q$ is a cubic root of unity, then $L$ is called a ${\mathbb Z}_3$-graded derivation.
\end{defin}

\subsection{Modules of ${\mathbb Z}_3$-graded algebras}\label{subsec2.2}

Since a general algebra do not need to have invertible elements, modules do not need always to have bases. In the ${\mathbb Z}_3$-graded case, there is an extra requirement of compatible degree.

\begin{defin} \label{def2.6}
Let ${\cal A}$ be a ${\mathbb Z}_3$-graded algebra and $M$ be a ${\mathbb Z}_3$-graded vector space. If there exists a mapping
${\cal A}\times M \longrightarrow M, \quad  (a,m) \mapsto am$
such that
$$p(am) = p(a) + p(m) \quad {\rm and} \quad  a(bm) = (ab) m $$
for all $a,b \in {\cal A}$ and all $m\in M$, then $M$ is called a left ${\mathbb Z}_3$-graded ${\cal A}$-module.
\end{defin}
A left ${\mathbb Z}_3$-graded ${\cal A}$-module $M$ over a ${\mathbb Z}_3$-graded algebra ${\cal A}$ is called free, if there exists a homogeneous basis$\{e_1,\ldots,e_k,e_{k+1},\ldots,e_{k+l},e_{k+l+1},\ldots,e_{k+l+n}\}$
for $M$, where $e_1,\ldots,e_k$ are elements of $M_0$, $e_{k+1},\ldots,e_{k+l}$ are elements of $M_1$ and $e_{k+l+1},\ldots,e_{k+l+n}$ are elements of $M_2$. So, every element $m\in M$ can uniquely be expressed as
$$m=\sum_i^{k+l+n} a^i e_i, \quad a^i \in {\cal A}.$$

For a ${\mathbb Z}_3$-graded algebra ${\cal A}$, there is a difference between left and right ${\cal A}$-modules. A left module can also be given the structure of a right module by defining the map
$$V\times {\cal A} \longrightarrow V, \quad  (v,a) \mapsto q^{p(v)p(a)} av.  $$
The set of ${\mathbb Z}_3$-graded derivations of ${\cal A}$ is an important example of a ${\mathbb Z}_3$-graded ${\cal A}$-module.

\subsection{${\mathbb Z}_3$-graded Hopf algebras}\label{subsec2.3}
The ${\mathbb Z}_2$-graded tensor product of two ${\mathbb Z}_3$-graded algebras ${\cal A}$ and ${\cal B}$ may be regarded as a ${\mathbb Z}_3$-graded algebra ${\cal A}\otimes{\cal B}$ with a product rule determined in \cite{Majid}
\begin{defin} \label{def2.7}
If ${\cal A}$ and ${\cal B}$ are two ${\mathbb Z}_3$-graded algebras, then the product rule in the ${\mathbb Z}_3$-graded algebra ${\cal A}\otimes {\cal B}$ is defined by
\begin{equation*}
(a_1\otimes b_1)(a_2\otimes b_2) = q^{p(b_1)p(a_2)} a_1a_2\otimes b_1b_2
\end{equation*}
where $a_i$'s and $b_i$'s are homogeneous elements in the algebras ${\cal A}$ and ${\cal B}$, respectively.
\end{defin}

\begin{defin} \label{def2.8}
A ${\mathbb Z}_3$-graded Hopf algebra is a ${\mathbb Z}_3$-graded vector space ${\cal A}$ over ${\mathbb K}$ with three linear map $\Delta$, $\epsilon$ and $S$ such that
\begin{align*}
(\Delta \otimes {\rm id}) \circ \Delta &= ({\rm id} \otimes \Delta) \circ \Delta, \\
m\circ(\epsilon\otimes{\rm id})\circ\Delta &= {\rm id} = m\circ({\rm id}\otimes\epsilon)\circ\Delta, \\
m\circ(S\otimes{\rm id})\circ\Delta &= \eta\circ\epsilon = m\circ({\rm id}\otimes S)\circ\Delta,
\end{align*}
together with $\Delta({\bf 1}) ={\bf 1}\otimes{\bf 1}$, $\epsilon({\bf 1})=1$, $S({\bf 1})={\bf 1}$, where $m$ is the product map, ${\rm id}$ is the identity map and $\eta:{\mathbb K}\longrightarrow{\cal A}$.
\end{defin}

\begin{remark} \label{rem2.1} 
The coproduct $\Delta$ is an algebra homomorphism from ${\cal A}$ to ${\cal A}\otimes{\cal A}$ which is multiplied in the way given in Definition 2.7 and the counit $\epsilon$ is an algebra homomorphism from ${\cal A}$ to ${\mathbb K}$.
\end{remark}

\begin{remark} \label{rem2.2} 
Any element of a ${\mathbb Z}_3$-graded Hopf algebra ${\cal A}$ is expressed as a product on the generators and its antipode (coinverse) $S$ is calculated with the property
\begin{equation*}
S(ab) = q^{p(a)p(b)} S(b)S(a), \qquad \forall a,b \in {\cal A}
\end{equation*}
in terms of antipode of the generators, where $q$ is a cubic root of unity.
\end{remark}

\subsection{${\mathbb Z}_3$-graded matrices}\label{subsec2.4}

If ${\cal A}$ is a ${\mathbb Z}_3$-graded algebra, ${\mathbb Z}_3$-graded matrices with entries in ${\cal A}$ define even homomorphisms of free
${\mathbb Z}_3$-graded ${\cal A}$-modules in terms of particular bases.

\begin{defin} \label{def2.9} 
An $n\times n$ matrix $T$ over a ${\mathbb Z}_3$-graded algebra ${\cal A}$ is a ${\mathbb Z}_3$-graded matrix whose entries are elements of ${\cal A}$ and which has the form
$T = T_0 + T_1 + T_2,$
where $T_0$, $T_1$ and $T_2$ are of grade 0, 1 and 2, respectively.
\end{defin}

\section{The ${\mathbb Z}_3$-graded quantum group $\widetilde{\rm GL}_q(2)$}\label{sec3}

In this section, we recall the ${\mathbb Z}_3$-graded quantum group $\widetilde{\rm GL}_q(2)$ from \cite{Celik4}. Let ${\cal O}(\widetilde{M}_q(2))$ be the complex associative algebra with the generators $a,\beta,\gamma,d$ satisfying the following relations:
\begin{align} \label{3.1}
a \beta &= \beta a, \qquad a \gamma = q \gamma a, \nonumber\\
\beta \gamma &= \gamma \beta, \qquad d\gamma = q^2 \gamma d, \\
\beta d &= d\beta, \qquad  ad = da + (q-1) \beta \gamma, \nonumber
\end{align}
where $q$ is a cubic root of unity. We assume that the generators $a$ and $d$ are of grade 0, the generators $\gamma$ and $\beta$ are of grade 1 and 2, respectively.

\begin{example} 
If we consider the generators of the algebra ${\cal O}(\widetilde{M}_q(2))$ as linear maps, then we can find the matrix representations of them. In fact,
it can be seen that there exists a ${\mathbb Z}_3$-graded representation $\rho:{\cal O}(\widetilde{M}_q(2))\to \mbox{M}(4,{\mathbb C})$ such that matrices
\begin{align} \label{3.2}
\rho(a) &= \begin{pmatrix} q & 0 & 0 & 0\\ 0 & q & 0 & 0\\ 0 & 0 & 1 & 0\\ 0 & 0 & 0 & 1 \end{pmatrix}, \quad
\rho(\beta)=\begin{pmatrix} 0 & q^2-q & 0 & 0\\ 0 & 0 & 0 & 0\\ 0 & 0 & 0 & q^2-q\\ 0 & 0 & 0 & 0\end{pmatrix}, \nonumber\\
\rho(\gamma) &= \begin{pmatrix} 0 & 0 & 1-q^2 & 0\\ 0 & 0 & 0 & 1-q^2 \\ 0 & 0 & 0 & 0\\ 0 & 0 & 0 & 0\end{pmatrix} , \quad
\rho(d)=\begin{pmatrix} 1 & 0 & 0 & q^2(q-1)^2\\ 0 & 1 & 0 & 0\\ 0 & 0 & q & 0\\ 0 & 0 & 0 & q \end{pmatrix}
\end{align}
representing the coordinate functions satisfy relations (\ref{3.1}), where $q$ is a cubic root of unity.
\end{example}

\noindent{\it Proof}.
Here we will show that the last relation in (\ref{3.1}) is realized, and it can be shown in a similar manner that other relations are realized. Since
\begin{equation*}
\rho_{ij}(fg) = q^{p(f)p(\rho_{ij}(g))} \sum_{k=1}^4 \rho_{ik}(f)\rho_{kj}(g), \qquad \forall f,g \in {\cal O}(\widetilde{M}_q(2))
\end{equation*}
we have $\rho_{ii}(ad) = q = \rho_{ii}(da)$ for each $i$, $\rho_{14}(ad) = (q-1)^2 = q\rho_{14}(da)$ and $\rho_{ij}(ad) = 0 = \rho_{ij}(da)$ for other $i$'s and $j$'s. Similarly, we find
\begin{equation*}
\rho_{14}(\beta\gamma) = q^2 \rho_{12}(\beta)\rho_{24}(\gamma) = q^2(q^2-q)(1-q^2)
\end{equation*}
and $\rho_{ij}(\beta\gamma) = 0$ for other $i$'s and $j$'s. Hence
\begin{align*}
\rho_{14}(ad-da-(q-1)\beta\gamma) &= (1-q^2)(q-1)^2 - (q-1)^2(1-q^2) = 0, \\
\rho_{ij}(ad-da-(q-1)\beta\gamma) &= q-q = 0
\end{align*}
for other $i$'s and $j$'s, so that we get
\begin{equation*}
\rho_{ij}(ad-da-(q-1)\beta\gamma) = 0, \quad \forall i,j
\end{equation*}
as expected. \hfill$\square$

\begin{example} 
Assume that $P\in\{q,q^2\}$ and $Q\in\{1,q^2\}$ where $q$ is a cubic root of unity. Then the matrices
\begin{align*}
L(a) &= \left(\begin{matrix} P & 0 & 0 & 0\\ 0 & P & 0 & 0\\ 0 & 0 & qP^{-1} & 0\\ 0 & 0 & 0 & P \end{matrix}\right), \quad
L(\beta) = \left(\begin{matrix} 0 & q(P-1) & 0 & 0\\ 0 & 0 & 0 & 0\\ 0 & 0 & 0 & 0\\ 0 & q^2(Q^{-1}-Q) & 0 & 0 \end{matrix}\right), \\
L(\gamma) &=\left(\begin{matrix} 0 & 0 & q^2P-P^{-1}) & 0\\ 0 & 0 & 0 & 0\\ 0 & 0 & 0 & 0\\ 0 & 0 & 0 & 0 \end{matrix}\right), \quad
L(d) = \left(\begin{matrix} 1 & 0 & 0 & P-1 \\ 0 & Q & 0 & 0\\ 0 & 0 & q & 0\\ 0 & 0 & 0 & qQ^{-1} \end{matrix}\right). \nonumber
\end{align*}
are a representation of $\widetilde{M}_q(2)$ with entries in ${\mathbb C}$, that is, they satisfy the relations (\ref{3.1}).
\end{example}

\begin{thm} \cite{Celik4} \label{theo3.1} 
There exists a unique bialgebra structure on the algebra ${\cal O}(\widetilde{M}_q(2))$ with the costructures
\begin{align} \label{3.3}
\Delta &:{\cal O}(\tilde{M}_q(2))\longrightarrow{\cal O}(\tilde{M}_q(2))\otimes{\cal O}(\tilde{M}_q(2)),\quad \Delta(t_{ij})=\sum_{k=1}^2 t_{ik}\otimes t_{kj}, \\
\epsilon &: {\cal O}(\tilde{M}_q(2))\longrightarrow {\mathbb C}, \quad \epsilon(t_{ij}) = \delta_{ij}
\end{align}
where $t_{11}=a$, $t_{12}=\beta$, $t_{21}=\gamma$, $t_{22}=d$.
\end{thm}
In addition, we have $\Delta({\bf 1})={\bf 1}\otimes {\bf 1}$ and $\epsilon({\bf 1}) = 1$.

The ${\mathbb Z}_3$-graded quantum determinant is defined in \cite{Celik4} by
\begin{equation} \label{3.5}
{\cal D}_q := ad - q \beta \gamma = da - \beta \gamma.
\end{equation}

\begin{remark} \label{rem3.2} 
The ${\mathbb Z}_3$-graded quantum determinant defined in (\ref{3.5}) is a group-like element belonging to the center of the algebra
${\cal O}(\widetilde{M}_q(2))$
\end{remark}

Using the quantum determinant ${\cal D}_q$ belonging to  the algebra ${\cal O}(\tilde{M}_q(2))$, we can define a Hopf algebra by adding the inverse
${\cal D}_q^{-1}$ to ${\cal O}(\tilde{M}_q(2))$. Let ${\cal O}(\widetilde{\rm GL}_q(2))$ be the quotient of the algebra ${\cal O}(\tilde{M}_q(2))$ by the two-sided ideal generated by the element $t{\cal D}_q-1$. For short we write
\begin{equation*}
{\cal O}(\widetilde{\rm GL}_q(2)) := {\cal O}(\tilde{M}_q(2))[t]/\langle t{\cal D}_q-1\rangle.
\end{equation*}
Then the algebra ${\cal O}(\widetilde{\rm GL}_q(2))$ is again a bialgebra.

\begin{thm} \cite{Celik4}\label{theo3.2} 
The bialgebra ${\cal O}(\widetilde{\rm GL}_q(2))$ is a ${\mathbb Z}_3$-graded Hopf algebra. The antipode $S$ of ${\cal O}(\widetilde{\rm GL}_q(2))$ is given by
\begin{equation} \label{3.6}
S(a) = d \,{\cal D}_q^{-1}, \quad S(\beta) = - \beta \,{\cal D}_q^{-1}, \quad S(\gamma) = -q\gamma \,{\cal D}_q^{-1}, \quad S(d) = a \,{\cal D}_q^{-1}.
\end{equation}
\end{thm}
In addition, we have $S({\bf 1}) = {\bf 1}$.

\begin{defin} \label{def3.2} 
The ${\mathbb Z}_3$-graded Hopf algebra ${\cal O}(\widetilde{\rm GL}_q(2))$ is called the coordinate algebra of the ${\mathbb Z}_3$-graded (quantum) group $\widetilde{\rm GL}_q(2)$.
\end{defin}

\begin{remark} \label{rem3.1} 
The cubes of the generators of the algebra ${\cal O}(\widetilde{\rm GL}_q(2))$ are central elements of ${\cal O}(\widetilde{\rm GL}_q(2))$.
\end{remark}

\section{Differential calculi on $\widetilde{\rm GL}_q(2)$}\label{sec4}

In this section, we will set up left-covariant differential calculi on the associative unital algebra ${\cal O}(\widetilde{\rm GL}_q(2))$ generated by four generators and quadratic and cubic relations. We need therefore to find commutation relations between the elements of the set
$\{a, \beta, \gamma, d, {\sf d}a, {\sf d}\beta, {\sf d}\gamma, {\sf d}d, {\sf d}^2a, {\sf d}^2\beta, {\sf d}^2\gamma, {\sf d}^2d\}$.

\subsection{${\mathbb Z}_3$-graded differential algebra} \label{sec4.1}

We know that, if ${\cal A}$ is an (associative) algebra then the central object of first order differential calculus is the exterior derivative
${\sf d}: {\cal A} \longrightarrow \Omega^1({\cal A})$ satisfying the Leibniz rule. The space of 1-forms $\Omega^1({\cal A})$ is generated as an
${\cal A}$-bimodule by the differentials of the generators of ${\cal A}$. It is furthermore required that the differentials form a basis of $\Omega^1({\cal A})$ as a left ${\cal A}$-module. To achieve this, we need to find commutation relations between the generators of ${\cal A}$ and their differentials which are consistent with the differential algebra structure.

\begin{defin} \label{def4.1}
Let ${\cal A}$ be a unital associative algebra, and $\Omega=\bigoplus_{i \in {\mathbb Z}_3} \Omega^{i}({\cal A})$ be a ${\mathbb Z}$-graded algebra with
$\Omega^0({\cal A})={\cal A}$. A pair $(\Omega,{\sf d})$ is called a ${\mathbb Z}_3$-graded differential calculus over ${\cal A}$ if ${\sf d}:\Omega\to\Omega$ is a linear map of degree one which satisfies
\renewcommand{\labelenumi}{(\arabic{enumi})}
\begin{enumerate}
\item ${\mathbb Z}_3$-graded Leibniz rule:
\begin{equation} \label{4.1}
{\sf d}(f\cdot g) = ({\sf d}f)\cdot g + q^{p(f)} f\cdot ({\sf d}g)
\end{equation}
for all homogeneous $f\in {\cal A}$ and all $g\in {\cal A}$, where $q$ is a primitive cubic root of unity,
\item ${\sf d}^3 = 0, \quad ({\sf d}^2\ne0)$.
\end{enumerate}
\end{defin}

In order to construct a differential calculus on ${\cal O}(\widetilde{\rm GL}_q(2))$, a linear operator {\sf d} which acts on the generators of the algebra
${\cal O}(\widetilde{\rm GL}_q(2))$ must be defined. For the definition, it is sufficient to define an action of {\sf d} on the generators and on their products as follows: for example, the operator {\sf d} is applied to $a$ produces a 1-form of degree 1, by definition. Similarly, application of {\sf d} to $\beta$ produces 1-form whose ${\mathbb Z}_3$-grade is 0. When the operator {\sf d} is applied twice by iteration to $a$ and $\beta$ it will produce a new entities which we will call a 2-form of degree 2 and 1, denoted by ${\sf d}^2a$ and ${\sf d}^2\beta$, respectively. Finally, we require that ${\sf d}^3=0$.

\subsection{${\mathbb Z}_3$-graded quantum de Rham complex of ${\cal O}(\widetilde{\rm GL}_q(2))$} \label{sec4.2}

The ${\mathbb Z}_3$-graded quantum de Rham complex of the algebra ${\cal O}(\widetilde{\rm GL}_q(2))$ consists of a first order differential calculus on
${\cal O}(\widetilde{\rm GL}_q(2))$ and its higher order extensions. Therefore, we need to find the relations between the generators of the algebra
${\cal O}(\widetilde{\rm GL}_q(2))$ and their first and second order differentials.

To obtain $q$-commutation relations between the generators of ${\cal O}(\widetilde{\rm GL}_q(2))$ and their first order differentials, we combine elements $a$, $\beta$, $\gamma$, $d$ and their differentials ${\sf d}a$, ${\sf d}\beta$, ${\sf d}\gamma$, ${\sf d}d$ which are considered as elements of a space
$\Omega^1 := \Omega^1({\cal O}(\widetilde{\rm GL}_q(2)))$ of 1-forms. We allow multiplication of the first order differentials in $\Omega^1$ from the left and from the right by the elements of ${\cal O}(\widetilde{\rm GL}_q(2))$, so that by the definition of the multiplication the resulting 1-form belongs to $\Omega^1$. This means that $\Omega^1$ is an ${\cal O}(\widetilde{\rm GL}_q(2))$-bimodule. The appearance of higher order differentials is a peculiar property of the ${\mathbb Z}_3$-graded calculus. The second order differentials ${\sf d}^2a$, ${\sf d}^2\beta$, ${\sf d}^2\gamma$, ${\sf d}^2d$ are considered as generators of a space $\Omega^2 := \Omega^2({\cal O}(\widetilde{\rm GL}_q(2)))$ of 2-forms which is actually also a bimodule.

\begin{defin} \label{def4.2}
A differential calculus $\Omega(\widetilde{\rm GL}_q(2))$ over the ${\mathbb Z}_3$-graded quantum group $\widetilde{\rm GL}_q(2)$ 
is called {\sf left covariant} with respect to ${\mathbb Z}_3$-graded Hopf algebra ${\cal O}(\widetilde{\rm GL}_q(2))$ if there exists a bimodule homomorphism
\begin{equation*}
\Delta_L:\Omega^1\longrightarrow {\cal O}(\widetilde{\rm GL}_q(2))\otimes\Omega^1
\end{equation*}
which is a left coaction of ${\cal O}(\widetilde{\rm GL}_q(2))$ on $\Omega^1$ such that
\begin{equation} \label{4.2}
\Delta_L(awb) = \Delta(a)\Delta_L(w)\Delta(b), \quad \forall a,b \in {\cal O}(\widetilde{\rm GL}_q(2)), \,\, w\in\Omega^1
\end{equation}
and
\begin{equation} \label{4.3}
\Delta_L\circ{\sf d} = (\tau'\otimes {\sf d})\circ \Delta,
\end{equation}
where $\tau':\Omega^1({\cal O}(\widetilde{\rm GL}_q(2)))\longrightarrow\Omega^1({\cal O}(\widetilde{\rm GL}_q(2)))$ is the linear map of degree zero which gives $\tau'(w)=q^{p(w)}w$ for all $w\in\Omega^1({\cal O}(\widetilde{\rm GL}_q(2)))$.
\end{defin}

Firstly, we will construct ${\mathbb Z}_3$-graded first order differential calculi on the algebra ${\cal O}(\widetilde{\rm GL}_q(2))$. Let the
${\cal O}(\widetilde{\rm GL}_q(2))$-bimodule $\Omega^1$ be generated as a free right ${\cal O}(\widetilde{\rm GL}_q(2))$-module by the differentials
${\sf d}a$, ${\sf d}\beta$, ${\sf d}\gamma$, ${\sf d}d$. In general the coordinates will not commute with their differentials. Therefore, one assumes that the possible commutation relations of the generators with their first order differentials are of the form
\begin{align} \label{4.4}
u_i\cdot {\sf d}v_j = \sum_{k,l=1}^4 B_{ij}^{kl} \, {\sf d}v_k\cdot u_l
\end{align}
where the constants $B_{ij}^{kl}$ are possibly depending on $q$. That is, let the left ${\cal O}(\widetilde{\rm GL}_q(2))$-module structure on $\Omega^1$ be completely defined by (\ref{4.4}).

\begin{remark} \label{rem4.1} 
As will be expected, there are eighty-six constants to be determined in the relations (\ref{4.4}) because the generators of the algebra
${\cal O}(\widetilde{\rm GL}_q(2))$ and their differentials are of graded. Therefore, the other constants contained therein should be considered to be automatically zero. For example,
\begin{equation*}
a\cdot {\sf d}\beta = B_{12}^{12} \, {\sf d}a\cdot\beta + B_{12}^{21} \, {\sf d}\beta\cdot a + B_{12}^{24} \, {\sf d}\beta\cdot d +
  B_{12}^{33} \, {\sf d}\gamma\cdot\gamma + B_{12}^{42} \, {\sf d}d\cdot\beta,
\end{equation*}
so that the constants $B_{12}^{11},B_{12}^{31},B_{12}^{41}$, {\it etc.} must be equal to zero.
\end{remark}

The existence of $\Delta_L$ depends on the commutation relations between $a$, $\beta$, $\gamma$, $d$ and their differentials. The following lemma
will help the proof of Theorem 4.5 below.
\begin{lemma} \label{lem4.1} 
Application of $\Delta_L$ to the relations $(\ref{4.4})$ leaves only ten constants in $(\ref{4.4})$.
\end{lemma}
\noindent{\it Proof}.
Using the identity in (\ref{4.3}), we can write
\begin{align} \label{4.5}
\Delta_L({\sf d}a) &= a\otimes {\sf d}a + q^2\beta\otimes {\sf d}\gamma, \quad
  \Delta_L({\sf d}\beta) = a\otimes {\sf d}\beta + q^2\beta\otimes {\sf d}d, \nonumber\\
\Delta_L({\sf d}\gamma) &= q\gamma\otimes {\sf d}a + d\otimes {\sf d}\gamma, \quad
  \Delta_L({\sf d}d) = q\gamma\otimes {\sf d}\beta +  d\otimes {\sf d}d.
\end{align}
Now we can apply the operator $\Delta_L$ from the left and right sides to the relations (\ref{4.4}) and obtain the following restrictions on some constants \begin{align} \label{4.6}
B_{13}^{13} &= B_{11}^{11} - q^2B_{13}^{31}, \quad B_{31}^{13} = B_{13}^{31}, \quad B_{31}^{31} = q^2(B_{11}^{11} - B_{13}^{31}), \quad
  B_{33}^{33} = q^2B_{11}^{11}, \nonumber\\
B_{41}^{14} &= q^2B_{23}^{32}, \quad B_{41}^{41} = (q^2-1)B_{21}^{21} - B_{23}^{23}, \quad B_{23}^{14} = q^2(B_{21}^{12} - B_{23}^{32}), \nonumber\\
B_{23}^{41} &= B_{21}^{21} - B_{23}^{23}, \quad B_{24}^{42} = B_{22}^{22} - B_{24}^{24}, \quad B_{32}^{23} = q^2B_{12}^{21}, \quad B_{32}^{32} = q^2B_{12}^{12},\\ B_{34}^{43} &= B_{12}^{21}, \quad B_{34}^{34} = q^2B_{12}^{12}, \quad B_{41}^{23} = B_{21}^{21} - B_{23}^{23}, \quad
  B_{41}^{32} = q^2B_{21}^{12} - B_{23}^{32}, \nonumber\\
B_{42}^{24} &= B_{22}^{22} - B_{24}^{24}, \quad B_{42}^{42} = (q^2-1)B_{22}^{22} + B_{24}^{24}, \,\,\, B_{43}^{34} = q^2B_{21}^{12}, \,\,\,
  B_{43}^{43} = q^2 B_{21}^{21}. \nonumber
\end{align}
The another constants in (\ref{4.4}) are equal to zero. \hfill$\square$

\begin{thm} \label{theo4.1} 
There exist left covariant ${\mathbb Z}_3$-graded first order differential calculi $\Omega(\widetilde{\rm GL}_q(2))$ over the algebra
${\cal O}(\widetilde{\rm GL}_q(2))$ with respect to itself such that $\{{\sf d}a,{\sf d}\beta,{\sf d}\gamma,{\sf d}d\}$ is a free right
${\cal O}(\widetilde{\rm GL}_q(2))$-module basis of $\Omega(\widetilde{\rm GL}_q(2))$. The bimodule structures for these calculi are determined by the relations
\begin{align} \label{4.7}
{\rm Case \, I}:\nonumber\\
a\cdot {\sf d}a &= P_1 \, {\sf d}a\cdot a, \quad
   a\cdot {\sf d}\beta = {\sf d}\beta\cdot a + (q^2P_2-1) \, {\sf d}a\cdot\beta, \nonumber\\
a\cdot {\sf d}\gamma &= q^2P_1^{-1} \, {\sf d}\gamma\cdot a + q(1-qP_1) \, {\sf d}a\cdot\gamma, \nonumber\\
a\cdot {\sf d}d &= {\sf d}d\cdot a + (q^2P_2-1) \, {\sf d}a\cdot d + (q-1) \, {\sf d}\beta\cdot\gamma, \nonumber\\
\beta\cdot {\sf d}a &= P_2 \, {\sf d}a\cdot \beta, \quad \beta\cdot {\sf d}\beta = P_2 \, {\sf d}\beta\cdot\beta, \quad
   \beta\cdot {\sf d}\gamma = P_2 \, {\sf d}\gamma\cdot\beta, \nonumber\\
\beta\cdot {\sf d}d &= P_2^{-1} \, {\sf d}d\cdot\beta + q^2(1-P_2){\sf d}\beta\cdot d, \\
\gamma\cdot {\sf d}a &= P_1^{-1} \, {\sf d}a\cdot\gamma + (q-P_1) \, {\sf d}\gamma\cdot a, \quad
   \gamma\cdot {\sf d}\gamma = q^2 P_1 \, {\sf d}\gamma\cdot\gamma, \nonumber\\
\gamma\cdot {\sf d}\beta &= q^2 \, {\sf d}\beta\cdot\gamma + q(P_2-q) \, {\sf d}\gamma\cdot\beta, \quad
   \gamma\cdot {\sf d}d = {\sf d}d\cdot\gamma + q(P_2-q) \, {\sf d}\gamma\cdot d, \nonumber\\
d\cdot {\sf d}a &= q^2P_2 \, {\sf d}a\cdot d + (q^2-1)P_2 \, {\sf d}\gamma\cdot\beta, \quad d\cdot {\sf d}\gamma = qP_2 \, {\sf d}\gamma\cdot d, \nonumber\\
d\cdot {\sf d}\beta &= P_2^{-1} \, {\sf d}\beta\cdot d + (q^2-P_2) \, {\sf d}d\cdot\beta, \quad
   d\cdot {\sf d}d = q^2 P_2 \, {\sf d}d\cdot d, \nonumber
\end{align}
\begin{align} \label{4.8}
{\rm Case \, II}:\nonumber\\
a\cdot {\sf d}a &= P_1 \, {\sf d}a\cdot a, \quad
  a\cdot {\sf d}\gamma = q^2P_1^{-1} \, {\sf d}\gamma\cdot a + q(1-qP_1) \, {\sf d}a\cdot\gamma, \nonumber\\
a\cdot {\sf d}\beta &= P_1 \, {\sf d}\beta\cdot a, \quad
  a\cdot {\sf d}d = P_1 \, {\sf d}d\cdot a + (q-1)P_1 \, {\sf d}\beta\cdot\gamma, \nonumber\\
\beta\cdot {\sf d}a &= q \, {\sf d}a\cdot \beta + q(P_1-1) \, {\sf d}\beta\cdot a, \quad
   \beta\cdot {\sf d}\beta = P_2 \, {\sf d}\beta\cdot\beta, \nonumber\\
\beta\cdot{\sf d}\gamma &= q \, {\sf d}\gamma\cdot\beta + q(P_1-1) \, {\sf d}\beta\cdot\gamma, \nonumber\\
\beta\cdot{\sf d}d &= P_2^{-1} \, {\sf d}d\cdot\beta + q^2(1-P_2) \, {\sf d}\beta\cdot d, \nonumber\\
\gamma\cdot{\sf d}a &= P_1^{-1} \, {\sf d}a\cdot\gamma + (q-P_1) \, {\sf d}\gamma\cdot a, \quad
   \gamma\cdot {\sf d}\beta = q^2P_1 \, {\sf d}\beta\cdot\gamma, \\
\gamma\cdot{\sf d}\gamma &= q^2P_1 \, {\sf d}\gamma\cdot\gamma, \quad
   \gamma\cdot{\sf d}d = P_1 \, {\sf d}d\cdot\gamma, \nonumber\\
d\cdot{\sf d}a &= {\sf d}a\cdot d + (1-q) \, {\sf d}\gamma\cdot\beta + (P_1-1) \, {\sf d}d\cdot a, \nonumber\\
d\cdot{\sf d}\beta &= P_2^{-1} \, {\sf d}\beta\cdot d + (q^2 - P_2) \, {\sf d}d\cdot\beta, \nonumber\\
d\cdot{\sf d}\gamma &= q^2 \, {\sf d}\gamma\cdot d + (P_1-1) \, {\sf d}d\cdot\gamma, \quad
   d\cdot {\sf d}d = q^2P_2 \, {\sf d}d\cdot d \nonumber
\end{align}
where $B_{11}^{11}=P_1\in\{q,q^2\}$ and $B_{22}^{22}=P_2\in\{1,q^2\}$ for both cases.
\end{thm}

\noindent{\it Proof}.
We apply the differential {\sf d} from the left to the relations (\ref{3.1}) and use relations associated with them from cross-commutation relations (\ref{4.4}) (with (\ref{4.6})). Later, we commute the differentials from the right to the left on both sides of the relations (\ref{3.1}). These enforce
\begin{align*}
{\rm Case \, I}: \,\, B_{12}^{21} = 1, \quad B_{12}^{12} = q^2P_2-1, \quad
{\rm Case \, II}: \,\, B_{12}^{21} = P_1, \quad B_{12}^{12} = 0,
\end{align*}
where in both cases $P_1\in\{q,q^2\}$, $P_2\in\{1,q^2\}$. \hfill$\square$

\begin{remark} \label{rem4.2} 
We will see from Theorem 4.14, which will be given later, that either $P_1=q^2$, $P_2=1$ or $P_1=q$, $P_2=q^2$.
\end{remark}

\begin{remark} \label{rem4.3} 
We know, from Theorem 4.5, that a free right ${\cal O}(\widetilde{\rm GL}_q(2))$-module basis of ${\cal O}(\widetilde{\rm GL}_q(2))$-bimodule
$\Omega(\widetilde{\rm GL}_q(2))$ is the set $\{{\sf d}a, {\sf d}\beta, {\sf d}\gamma, {\sf d}d\}$ and the relations (\ref{4.7}) (and (\ref{4.8})) hold. We now consider left module structure of ${\cal O}(\widetilde{\rm GL}_q(2))$-bimodule $\Omega(\widetilde{\rm GL}_q(2))$. The left product
${\sf d}v \mapsto u\cdot{\sf d}v$ is an endomorphism of the right module $\Omega(\widetilde{\rm GL}_q(2))$. Ring of all endomorphisms of any free module of rank 4 is isomorphic to the ring of all $4\times4$-matrices. Since $\{{\sf d}a, {\sf d}\beta, {\sf d}\gamma, {\sf d}d\}$ is the homogeneous basis of
$\Omega(\widetilde{\rm GL}_q(2))$, there exists a map
$\sigma:{\cal O}(\widetilde{\rm GL}_q(2))\to M_4({\cal O}(\widetilde{\rm GL}_q(2)))$ defined by the formulas
\begin{equation} \label{4.9}
f\cdot{\sf d}g_j = \sum_{i=1}^4 q^{p(f)p({\sf d}g_j)} \,{\sf d}g_i\cdot\sigma_{ij}(f)
\end{equation}
for all $f\in {\cal O}(\widetilde{\rm GL}_q(2))$ and $g_1=a$, $g_2=\beta$, $g_3=\gamma$ and $g_4=d$. Indeed, one can see that the relations (\ref{4.9}) equivalent to the relations (\ref{4.7}), where
\begin{align} \label{4.10}
\sigma(a) &= \left(\begin{matrix} P_1a & (q^2P_2-1)\beta & q(1-qP_1)\gamma & (q^2P_2-1)d \\ 0 & a & 0 & (q-1)\gamma \\ 0 & 0 & q^2P_1^{-1}a & 0 \\ 0 & 0 & 0 & a\end{matrix}\right), \nonumber\\
\sigma(\beta) &= \left(\begin{matrix} qP_2\beta & 0 & 0 & 0 \\ 0 & P_2\beta & 0 & (1-P_2)d \\ 0 & 0 & q^2P_2\beta & 0 \\ 0 & 0 & 0 & qP_2^{-1}\beta \end{matrix}\right), \nonumber\\
\sigma(\gamma) &= \left(\begin{matrix} q^2P_1^{-1} \gamma & 0 & 0 & 0 \\ 0 & q^2\gamma & 0 & 0 \\ (1-q^2P_1)a & q(P_2-q) \beta & P_1\gamma & (P_2-q) d \\ 0 & 0 & 0 & q^2\gamma \end{matrix}\right), \\
\sigma(d) &= \left(\begin{matrix} q^2P_2 d & 0 & 0 & 0 \\ 0 & P_2^{-1}d & 0 & 0 \\ (q^2-1)P_2\beta & 0 & qP_2d & 0 \\ 0 & (q^2-P_2)\beta & 0 & q^2P_2d \end{matrix}\right) \nonumber
\end{align}
with $P_1\in\{q,q^2\}$ and $P_2\in\{1,q^2\}$.
\end{remark}

The formulae in (\ref{4.9}) and the matrices in (\ref{4.10}) are written for the Case I. Similar formulae and matrices can be repeated for the Case II.

\begin{thm} \label{theo4.2} 
The map $\sigma$ is ${\mathbb C}$-linear ${\mathbb Z}_3$-graded homomorphism such that
\begin{equation} \label{4.11}
\sigma_{ij}(fg) = \sum_{k=1}^4 q^{p(f)[p({\sf d}f_k)-p({\sf d}f_j)]} \,\sigma_{ik}(f)\sigma_{kj}(g), \quad \forall f,g\in {\cal O}(\widetilde{\rm GL}_q(2)).
\end{equation}
\end{thm}

\noindent{\it Proof}.
According to (\ref{4.9}), we can write
\begin{align*}
(uv)\cdot{\sf d}f_j = \sum_{i=1}^4 q^{p(uv)p({\sf d}f_j)} \,{\sf d}f_i\cdot\sigma_{ij}(uv), \quad \forall u,v\in {\cal O}(\widetilde{\rm GL}_q(2)).
\end{align*}
In the other hand, using the identity $(uv)\cdot {\sf d}f_i = u(v\cdot {\sf d}f_i)$, one obtains
\begin{align*}
(uv)\cdot{\sf d}f_j
&= q^{p(v)p({\sf d}f_j)} \, u\cdot \sum_{i=1}^4 {\sf d}f_i\cdot \sigma_{ij}(v)\\
&= q^{p(v)p({\sf d}f_j)} \, \sum_{i=1}^4 \left(q^{p(u)p({\sf d}f_i)} \, \sum_{k=1}^4 {\sf d}f_k \cdot \sigma_{ki}(u) \right) \sigma_{ij}(v) \\
&= \sum_{i=1}^4\sum_{k=1}^4 q^{p(v)p({\sf d}f_j) + p(u)p({\sf d}f_i)} \, {\sf d}f_k \cdot \sigma_{ki}(u) \sigma_{ij}(v).
\end{align*}
Interchanging $i$ and $k$ in the last expression gives the equality in (\ref{4.11}). \hfill$\square$

\begin{remark} \label{rem4.4} 
It is easy to see that the relations (\ref{3.1}) are preserved under the action of the operator $\sigma$.
\end{remark}

\begin{remark} \label{rem4.5} 
We can also define a map  $\tau:{\cal O}(\widetilde{\rm GL}_q(2))\to M_4({\cal O}(\widetilde{\rm GL}_q(2)))$ by the formulas
\begin{equation} \label{4.12}
{\sf d}g_i \cdot f= \sum_{j=1}^4 q^{p(f)p({\sf d}g_j)} \,\tau_{ij}(f)\cdot{\sf d}g_j
\end{equation}
for all $f\in {\cal O}(\widetilde{\rm GL}_q(2))$. For the Case I, we have
\begin{align} \label{4.13}
\tau(a) &= \left(\begin{matrix} P_1^{-1}a & 0 & 0 & 0\\ (P_2^{-1}-q^2)\beta & a & 0 & 0\\ (1-q^2P_1^{-1})\gamma & 0 & qP_1a & 0\\
(qP_2^{-1}-1)d & (q-q^2)\gamma & 0 & a \end{matrix}\right), \nonumber\\
\tau(\beta) &=\left(\begin{matrix} qP_2^{-1}\beta & 0 & 0 & 0\\ 0 & P_2^{-1}\beta & 0 & 0\\ 0 & 0 & q^2P_2^{-1}\beta & 0\\ 0 & q^2(1-P_2^{-1})d & 0 & qP_2\beta \end{matrix}\right), \nonumber\\
\tau(\gamma) &= \left(\begin{matrix} q^2P_1 \gamma & 0 & q^2(1-qP_1^{-1})a & 0 \\ 0 & q\gamma & (qP_2^{-1}-1)\beta & 0\\ 0 & 0 & q^2P_1^{-1}\gamma & 0\\ 0 & 0 & q(qP_2^{-1}-1)d & q^2\gamma \end{matrix}\right), \\
\tau(d) &= \left(\begin{matrix} qP_2^{-1}d & 0 & (q-1)P_2^{-1}\beta & 0\\ 0 & P_2d & 0 & (1-q^2P_2^{-1})\beta\\ 0 & 0 & q^2P_2^{-1}d & 0\\ 0 & 0 & 0 & qP_2^{-1}d \end{matrix}\right) \nonumber
\end{align}
where $P_1\in\{q,q^2\}$ and $P_2\in\{1,q^2\}$.
\end{remark}

\begin{thm} \label{theo4.3} 
The map $\tau$ is ${\mathbb C}$-linear ${\mathbb Z}_3$-graded homomorphism such that
\begin{equation} \label{4.14}
\tau_{ij}(fg) = \sum_{k=1}^4 q^{p(f)[p({\sf d}f_k)-p({\sf d}f_j)]} \,\tau_{ik}(f)\tau_{kj}(g), \quad \forall f,g\in {\cal O}(\widetilde{\rm GL}_q(2)).
\end{equation}
\end{thm}

\noindent{\it Proof}.
From (\ref{4.12}), we can write
\begin{align*}
{\sf d}f_i\cdot (uv) = \sum_{i=1}^4 q^{p(uv)p({\sf d}f_j)} \,{\sf d}f_i\cdot\sigma_{ij}(uv), \quad \forall u,v\in {\cal O}(\widetilde{\rm GL}_q(2)).
\end{align*}
In the other hand,
\begin{align*}
{\sf d}f_i\cdot (uv)
&= \left(\sum_{j=1}^4 q^{p(u)p({\sf d}f_j)} \, \tau_{ij}(u)\cdot {\sf d}f_j\right)\cdot v \\
&= \sum_{j=1}^4 q^{p(u)p({\sf d}f_j)} \, \tau_{ij}(u)\cdot \left(\sum_{k=1}^4 q^{p(v)p({\sf d}f_k)} \, \tau_{jk}(v) \cdot  {\sf d}f_k\right) \\
&= \sum_{j=1}^4\sum_{k=1}^4 q^{p(u)p({\sf d}f_k) + p(v)p({\sf d}f_j)} \, \tau_{ik}(u)\tau_{kj}(u)\cdot {\sf d}f_j
\end{align*}
so that we conclude the equality in (\ref{4.14}). \hfill$\square$

\begin{remark} \label{rem4.6} 
It is easy to see from Theorem 4.11 that the relations (\ref{3.1}) are preserved under the action of the operator $\tau$.
\end{remark}

\begin{corollary} \label{cor4.1} 
The differentiations of the cubes of the generators of the algebra ${\cal O}(\widetilde{\rm GL}_q(2))$ are all zero.
\end{corollary}

\noindent{\it Proof}.
Let us write $f\cdot{\sf d}f=Q\, {\sf d}f\cdot f$ where $f$ is an arbitrary element of ${\cal O}(\widetilde{\rm GL}_q(2))$. Then, using (\ref{4.1}),
we can write
\begin{align*}
{\sf d}(f^3)
&= {\sf d}f\cdot f^2 + q^{p(f)} f\cdot ({\sf d}f\cdot f+ q^{p(f)} f\cdot {\sf d}f) \\
&= (1 + q^{p(f)}Q + q^{2p(f)}Q^2) \, {\sf d}f\cdot f^2.
\end{align*}
Since
\begin{equation*}
1 + q^{p(f)}Q + q^{2p(f)}Q^2 = 0
\end{equation*}
for each $f\in{\cal O}(\widetilde{\rm GL}_q(2))$, we say that ${\sf d}(f^3)$ should be zero. \hfill$\square$

\begin{thm} \label{theo4.4} 
The commutation relations between the differentials are as follows
\begin{align} \label{4.15}
{\sf d}a\cdot{\sf d}\beta &= q^2 \, {\sf d}\beta\cdot{\sf d}a, \quad
  {\sf d}a\cdot{\sf d}\gamma = q^2 \, {\sf d}\gamma\cdot{\sf d}a, \nonumber\\
{\sf d}d\cdot{\sf d}\beta &= q^2 \, {\sf d}\beta\cdot{\sf d}d, \quad
  {\sf d}d\cdot{\sf d}\gamma = {\sf d}\gamma\cdot{\sf d}d, \\
{\sf d}\beta\cdot{\sf d}\gamma &= q^2 \, {\sf d}\gamma\cdot{\sf d}\beta, \quad
  {\sf d}a\cdot{\sf d}d = {\sf d}d\cdot{\sf d}a + (q^2-q) \, {\sf d}\gamma\cdot{\sf d}\beta.\nonumber
\end{align}
\end{thm}

\noindent{\it Proof}.
Inspired by the relations (\ref{3.1}), we can assume that the possible relations between the differentials are of the form
\begin{align} \label{4.16}
{\sf d}a\cdot{\sf d}\beta &= X_1 \, {\sf d}\beta\cdot{\sf d}a, \quad
  {\sf d}a\cdot{\sf d}\gamma = X_2 \, {\sf d}\gamma\cdot{\sf d}a, \nonumber\\
{\sf d}d\cdot{\sf d}\beta &= X_3 \, {\sf d}\beta\cdot{\sf d}d, \quad
  {\sf d}d\cdot{\sf d}\gamma = X_4 \, {\sf d}\gamma\cdot{\sf d}d, \\
{\sf d}a\cdot{\sf d}d &= X_5 \, {\sf d}d\cdot{\sf d}a + X_6 \, {\sf d}\gamma\cdot{\sf d}\beta, \quad
{\sf d}\beta\cdot{\sf d}\gamma = X_7 \, {\sf d}\gamma\cdot{\sf d}\beta + X_8 \, {\sf d}d\cdot{\sf d}a.\nonumber
\end{align}
If we now apply the operator $\Delta_L$ given in (\ref{4.5}) from the left and the right sides to the relations (\ref{4.16}), then we have
$X_2 = q^2$, $X_3 = q^2$, $X_4 = q^2X_1^{-1}$, $X_6 = q^2 - X_1^{-1}X_5$, $X_7 = -q^2X_1^{-1}(1+X_1^{-1}X_5)$, $X_8 = q^2X_1^{-1}X_5-1$. On the other hand, we check that all ambiguities of cubic monomials are resolvable with respect to any of the orders
${\sf d}a<{\sf d}d<{\sf d}\beta<{\sf d}\gamma<a<d<\beta<\gamma$ or $a<d<\beta<\gamma<{\sf d}a<{\sf d}d<{\sf d}\beta<{\sf d}\gamma$.
For example, $a\cdot({\sf d}\beta\cdot{\sf d}a)$ must be equal to $(a\cdot{\sf d}\beta)\cdot{\sf d}a$, {\it etc.} When the necessary processing is done, we see that $B_{12}^{21}=P_1=P_2X_1$ for the Case I and $B_{12}^{12}=q^2P_1X_1^{-1}-1$ for the Case II. Consequently, we have
$X_1 = P_1P_2^{-1}$; $X_4 = q^2X_1^{-1}$, $X_5 = qX_1$, $X_6 = q^2 - q$, $X_7 = qX_1^{-1}$, $X_8 = 0$ if $P_1=q^2$, $P_2=1$ or $P_1=q$, $P_2=q^2$. In both cases we have the relations (\ref{4.15}). \hfill$\square$

Now we want to find the relations between the generators of  ${\cal O}(\widetilde{\rm GL}_q(2))$, the first and second order differentials of the generators for the Case I. (Similar results can be repeated for the Case II). To do this, we first take the differentials of both sides of the equations in (\ref{4.7}).
These processes give relations between the generators and their second order differentials, including the first order differentials of the generators: for example,
\begin{equation*}
a\cdot{\sf d}^2\beta = {\sf d}^2\beta\cdot a + (q^2P_2-1) \, {\sf d}^2a\cdot\beta + (P_2-1) \, {\sf d}a\cdot{\sf d}\beta.
\end{equation*}
Such relations are not homogeneous in the sense that the commutation relations between the generators and the second order differentials include the first order differentials as well. In other words, the above equation indicates also that, e.g. ${\sf d}a\cdot{\sf d}\beta$ etc., seem to be in $\Omega^2$. However,
${\sf d}a\cdot{\sf d}\beta$ etc., depend on both ${\sf d}^2a\cdot\beta$ and ${\sf d}^2\beta\cdot a$ etc., and it will be discussed in Lemma 4.17.
In particular, we can assume that there exist commutation relations between the generators of ${\cal O}(\widetilde{\rm GL}_q(2))$ and their second order differentials as follows
\begin{align} \label{4.18}
u_i\cdot{\sf d}^2v_j = \sum C_{ij}^{kl} \, {\sf d}^2v_k\cdot u_l
\end{align}
where the constants $C_{ij}^{kl}$ are possibly depend on $q$. This request allows us to introduce a left ${\cal O}(\widetilde{\rm GL}_q(2))$-module, therefore a bimodule, structure on a right free module generated by the second order differentials of the generators.

\begin{lemma} \label{lem4.2} 
Let $P_2\in\{1,q^2\}$. Commutation relations of the first order differentials with the second order differentials are as follows
\begin{align}\label{4.19}
{\sf d}a\cdot{\sf d}^2a & = q \, {\sf d}^2a\cdot{\sf d}a, \quad {\sf d}a\cdot{\sf d}^2\gamma = {\sf d}^2\gamma\cdot{\sf d}a, \nonumber\\
{\sf d}a\cdot{\sf d}^2\beta &= qP_2 \, {\sf d}^2\beta\cdot{\sf d}a + qP_2(1-qP_2) \, {\sf d}^2a\cdot{\sf d}\beta, \nonumber\\
{\sf d}a\cdot{\sf d}^2d &= q^2P_2 \, {\sf d}^2d\cdot{\sf d}a + qP_2(1-qP_2) \, {\sf d}^2a\cdot{\sf d}d + P_2(q^2-q) \, {\sf d}^2\beta\cdot{\sf d}\gamma, \nonumber\\
{\sf d}\beta\cdot{\sf d}^2a &= P_2^{-1} \, {\sf d}^2a\cdot{\sf d}\beta + qP_2(P_2-1) \, {\sf d}^2\beta\cdot{\sf d}a, \quad
   {\sf d}\beta\cdot{\sf d}^2\beta = q^2 \, {\sf d}^2\beta\cdot{\sf d}\beta, \nonumber\\
{\sf d}\beta\cdot{\sf d}^2\gamma &= qP_2^{-1} \, {\sf d}^2\gamma\cdot{\sf d}\beta + qP_2(P_2-1) \, {\sf d}^2\beta\cdot{\sf d}\gamma, \quad
   {\sf d}\beta\cdot{\sf d}^2d = {\sf d}^2d\cdot{\sf d}\beta, \nonumber\\
{\sf d}\gamma\cdot{\sf d}^2a &= q^2 \, {\sf d}^2a\cdot{\sf d}\gamma + (1-q) \, {\sf d}^2\gamma\cdot{\sf d}a, \quad
   {\sf d}\gamma\cdot{\sf d}^2\gamma = {\sf d}^2\gamma\cdot{\sf d}\gamma, \nonumber\\
{\sf d}\gamma\cdot{\sf d}^2\beta &= q^2P_2 \, {\sf d}^2\beta\cdot{\sf d}\gamma + P_2(1-qP_2) \, {\sf d}^2\gamma\cdot{\sf d}\beta, \\
{\sf d}\gamma\cdot{\sf d}^2d & = qP_2 \, {\sf d}^2d\cdot{\sf d}\gamma + P_2(1-qP_2) \, {\sf d}^2\gamma\cdot{\sf d}d, \nonumber\\
{\sf d}d\cdot{\sf d}^2a &= qP_2^{-1} \, {\sf d}^2a\cdot{\sf d}d + P_2^{-1}(q^2-1) \, {\sf d}^2\gamma\cdot{\sf d}\beta + P_2(P_2-1) \, {\sf d}^2d\cdot{\sf d}a, \nonumber\\
{\sf d}d\cdot{\sf d}^2\beta &= q \, {\sf d}^2\beta\cdot{\sf d}d + (q-q^2) \, {\sf d}^2d\cdot{\sf d}\beta, \nonumber\\
{\sf d}d\cdot{\sf d}^2\gamma &= qP_2^{-1} \, {\sf d}^2\gamma\cdot{\sf d}d + P_2(P_2-1) \, {\sf d}^2d\cdot{\sf d}\gamma, \quad
   {\sf d}d\cdot{\sf d}^2d = q \, {\sf d}^2d\cdot{\sf d}d. \nonumber
\end{align}
\end{lemma}

\noindent{\it Proof}.
Twice application of the differential {\sf d} on the left and right sides of the relations in (\ref{4.7}) gives rise the relations in (\ref{4.19}). \hfill$\square$

\begin{corollary} \label{cor4.2} 
The differentiations of $({\sf d}f)^3$ for all $f\in{\cal O}(\widetilde{\rm GL}_q(2))$ are all zero.
\end{corollary}

\noindent{\it Proof}.
If $f$ is an arbitrary element of ${\cal O}(\widetilde{\rm GL}_q(2))$, using (\ref{4.1}) we obtain
\begin{align*}
{\sf d}(({\sf d}f)^3)
&= {\sf d}^2f\cdot ({\sf d}f)^2 + q^{p({\sf d}f)} {\sf d}f\cdot ({\sf d}^2f\cdot {\sf d}f + q^{p({\sf d}f)} {\sf d}f\cdot {\sf d}^2f) \\
&= (1 + q^{p({\sf d}f)-p({\sf d}^2f)} + q^{2[p({\sf d}f)-p({\sf d}^2f)]}) \, {\sf d}^2f\cdot ({\sf d}f)^2 = 0,
\end{align*}
since
\begin{equation*}
1 + q^{\alpha} + q^{2\alpha} = 0
\end{equation*}
where $\alpha=p({\sf d}f)-p({\sf d}^2f)$ for all $f\in{\cal O}(\widetilde{\rm GL}_q(2))$. \hfill$\square$

\newpage
\begin{lemma} \label{lem4.3} 
The products of the first order differentials are related with the second order differentials by the following equations
\begin{align}\label{4.20}
(qP_1-1) \, {\sf d}a\cdot{\sf d}a & = (q^2-P_1) \, {\sf d}^2a\cdot a, \quad (P_2-1) \, {\sf d}\beta\cdot{\sf d}\beta = (1-P_2) \, {\sf d}^2\beta\cdot\beta, \nonumber\\
(P_2-1) \, {\sf d}a\cdot{\sf d}\beta &= P_2(P_2-1) \, (q \, {\sf d}^2\beta\cdot a + P_2 \, {\sf d}^2a\cdot\beta), \nonumber\\
(q-P_1^{-1}) \, {\sf d}a\cdot{\sf d}\gamma &= (qP_1-1) \, ({\sf d}^2\gamma\cdot a + {\sf d}^2a\cdot\gamma), \nonumber\\
(P_2-1) \, {\sf d}a\cdot{\sf d}d &= (P_2-1) \, [q \, {\sf d}^2d\cdot a + P_2 \, {\sf d}^2a\cdot d +(q^2-q) \, {\sf d}^2\beta\cdot\gamma], \nonumber\\
(P_2-1) \, {\sf d}\beta\cdot{\sf d}\gamma &= P_2(P_2-1) \, (P_2 \, {\sf d}^2\gamma\cdot\beta + q^2 \, {\sf d}^2\beta\cdot\gamma), \\
(P_2-1) \, {\sf d}\beta\cdot{\sf d}d &= (P_2-1) \, (q \, {\sf d}^2d\cdot\beta + q^2 \, {\sf d}^2\beta\cdot d), \nonumber\\
(P_2-1) \, {\sf d}\gamma\cdot{\sf d}d &= P_2(P_2-1) \, (q \, {\sf d}^2d\cdot\gamma + q^2P_2 \, {\sf d}^2\gamma\cdot d), \nonumber\\
(P_1-q^2) \, {\sf d}\gamma\cdot{\sf d}\gamma & = q(q^2-P_1) \, {\sf d}^2\gamma\cdot\gamma, \quad
(P_2-1) \, {\sf d}d\cdot{\sf d}d = q^2(1-P_2) \, {\sf d}^2d\cdot d, \nonumber
\end{align}
where $P_1\in\{q,q^2\}$ and $P_2\in\{1,q^2\}$.
\end{lemma}

\noindent{\it Proof}.
When we take the differentials of the relations (\ref{4.18}) and compare them with (\ref{4.19}), we find both the constants $C_{ij}^{kl}$ and the relations (\ref{4.20}). \hfill$\square$

Note that, if $P_1=q^2$ and $P_2=1$ in the equalities in Lemma 4.17, then there exist no relations between the products of the first order differentials with
the products of the second order differentials with the generators. In this case, relations (\ref{4.21}) below can be easily obtained from the relations (\ref{4.7}).

\begin{corollary} \label{cor4.3} 
The relations between the elemens of ${\cal O}(\widetilde{\rm GL}_q(2))$ and their second order differentials are as follows
\begin{align} \label{4.21}
a\cdot {\sf d}^2a &= q^2 \, {\sf d}a\cdot a, \quad a\cdot{\sf d}^2\gamma = {\sf d}^2\gamma\cdot a, \nonumber\\
a\cdot{\sf d}^2\beta &= P_2 \, {\sf d}^2\beta\cdot a + P_2(q^2-P_2) \, {\sf d}^2a\cdot\beta, \nonumber\\
a\cdot{\sf d}^2d &= P_2 \, {\sf d}^2d\cdot a + P_2(q^2-P_2) \, {\sf d}^2a\cdot d + P_2(q-1) \, {\sf d}^2\beta\cdot\gamma, \nonumber\\
\beta\cdot{\sf d}^2a &= qP_2^2 \, {\sf d}^2a\cdot\beta + P_2(P_2-1) \, {\sf d}^2\beta\cdot a, \quad
   \beta\cdot{\sf d}^2\beta = q \, {\sf d}^2\beta\cdot\beta, \nonumber\\
\beta\cdot{\sf d}^2\gamma &= qP_2^{-1} \, {\sf d}^2\gamma\cdot\beta + P_2(P_2-1) \, {\sf d}^2\beta\cdot\gamma, \quad
   \beta\cdot{\sf d}^2d = q \, {\sf d}^2d\cdot\beta, \\
\gamma\cdot{\sf d}^2a &= {\sf d}^2a\cdot\gamma + (1-q) \, {\sf d}^2\gamma\cdot a, \quad
   \gamma\cdot{\sf d}^2\gamma = {\sf d}^2\gamma\cdot\gamma, \nonumber\\
\gamma\cdot{\sf d}^2\beta &= qP_2 \, {\sf d}^2\beta\cdot\gamma + P_2(1-qP_2) \, {\sf d}^2\gamma\cdot\beta, \nonumber\\
\gamma\cdot{\sf d}^2d &= q^2P_2 \, {\sf d}^2d\cdot\gamma + P_2(1-qP_2) \, {\sf d}^2\gamma\cdot d, \nonumber\\
d\cdot{\sf d}^2a &= q^2P_2^{-1} \, {\sf d}^2a\cdot d + P_2^{-1}(q^2-1) \, {\sf d}^2\gamma\cdot\beta + P_2(P_2-1) \, {\sf d}^2d\cdot a, \nonumber\\
d\cdot{\sf d}^2\beta &= {\sf d}^2\beta\cdot d + (q^2-1) \, {\sf d}^2d\cdot\beta, \quad
d\cdot{\sf d}^2d = q^2 \, {\sf d}d^2\cdot d, \nonumber\\
d\cdot{\sf d}^2\gamma &= qP_2^{-1} \, {\sf d}^2\gamma\cdot d + qP_2(P_2-1) \, {\sf d}^2d\cdot\gamma, \nonumber
\end{align}
where $P_2\in\{1,q^2\}$.
\end{corollary}

\begin{corollary} \label{cor4.4} 
The commutation relation between the second order differentials reads as follows
\begin{align} \label{4.22}
{\sf d}^2a\cdot{\sf d}^2\beta &= q \, {\sf d}^2\beta\cdot{\sf d}^2a, \quad
   {\sf d}^2\beta\cdot{\sf d}^2\gamma = q \, {\sf d}^2\gamma\cdot{\sf d}^2\beta, \nonumber\\
{\sf d}^2a\cdot{\sf d}^2\gamma &= {\sf d}^2\gamma\cdot{\sf d}^2a, \quad
   {\sf d}^2\beta\cdot{\sf d}^2d = q^2 \, {\sf d}^2d\cdot{\sf d}^2\beta, \\
{\sf d}^2\gamma\cdot{\sf d}^2d &= q^2 \, {\sf d}^2d\cdot{\sf d}^2\gamma, \quad
   {\sf d}^2a\cdot{\sf d}^2d = {\sf d}^2d\cdot{\sf d}^2a + (q^2-q) \, {\sf d}^2\beta\cdot{\sf d}^2\gamma. \nonumber
\end{align}
\end{corollary}

For the Case I, the relations (\ref{3.1}), (\ref{4.7}), (\ref{4.15}), (\ref{4.19}), (\ref{4.21}) and (\ref{4.22}) define a universal ${\mathbb Z}_3$-graded quantum de Rham complex of $\widetilde{\rm GL}_q(2)$. The same can be repeated for the Case II.

\subsection{The ${\mathbb Z}_3$-graded Maurer-Cartan forms} \label{sec4.3}

The left invariant Maurer-Cartan forms play a special role in the classical geometry of the Lie groups. To find their quantum analogues we first have to need convert the left action of group on forms to a left coaction. This situation is described by the following definition:

\begin{defin} \label{def4.3} 
Let $\Delta_L$ be a linear map from a bimodule $\Omega$ over ${\cal O}(\widetilde{\rm GL}_q(2))$ to ${\cal O}(\widetilde{\rm GL}_q(2))\otimes\Omega$. If, for all $f,g \in {\cal O}(\widetilde{\rm GL}_q(2))$ and $w\in\Omega$
\begin{align*}
\Delta_L(fwg) &= \Delta(f)\Delta_L(w)\Delta(g), \quad (\epsilon\otimes{\rm id})\Delta_L(w)=w \quad {\rm and}\\
(\Delta\otimes{\rm id})\circ\Delta_L &= ({\rm id}\otimes\Delta_L)\circ\Delta_L
\end{align*}
we say that $(\Omega,\Delta_L)$ is a left-covariant bimodule.
\end{defin}

We now introduce the left-invariant Maurer-Cartan forms as
\begin{equation} \label{4.23}
w_{ij}=:w(t_{ij}) = m\circ(S'\otimes{\sf d})\Delta(t_{ij}) \quad {\rm or} \quad w_{ij} = \sum_{k=1}^2 \,q^{p(t_{ik})} \, S(t_{ik})\cdot{\sf d}t_{ik}
\end{equation}
where $S'(t_{ij})=q^{p(t_{ij})} \, S(t_{ij})$. Explicitly,
\begin{align} \label{4.24}
w_1 &= {\cal D}_q^{-1} \, (d\cdot{\sf d}a - q^2\beta\cdot{\sf d}\gamma), \quad w_3 = {\cal D}_q^{-1} \, (a\cdot{\sf d}\gamma - q^2\gamma\cdot{\sf d}a), \nonumber\\
w_2 &= {\cal D}_q^{-1} \, (d\cdot{\sf d}\beta - q^2\beta\cdot{\sf d}d),  \quad w_4 = {\cal D}_q^{-1} \, (a\cdot{\sf d}d - q^2\gamma\cdot{\sf d}\beta),
\end{align}
where $w_{11}=w_1$, $w_{12}=w_2$, $w_{21}=w_3$ and $w_{22}=w_4$.

\begin{defin} \label{def4.4} 
A form $w$ is left-coinvariant if $\Delta_L(w) = {\bf 1}\otimes w$.
\end{defin}

\begin{thm} \label{theo4.5} 
The forms defined by in (\ref{4.24}) are coinvariant under the left-coaction.
\end{thm}

\noindent{\it Proof}.
As a sample let us check that the form $w_1$ is indeed left-coinvariant:
\begin{align*}
\Delta_L(w_1)
&= \Delta_L({\cal D}_q^{-1} \, (d\cdot{\sf d}a - q^2\beta\cdot{\sf d}\gamma)), \\
&= {\cal D}_q^{-1} \left[ (\gamma\otimes\beta + d\otimes d)(a\otimes{\sf d}a + q^2\beta\otimes{\sf d}\gamma)  \right.\\
   &\quad -\left. q^2 (a\otimes\beta + \beta\otimes d)(q\gamma\otimes{\sf d}a + d\otimes{\sf d}\gamma)\right]\\
&= {\cal D}_q^{-1} \left[\gamma a\otimes\beta\cdot{\sf d}a + \gamma\beta\otimes\beta\cdot{\sf d}\gamma + da\otimes d\cdot{\sf d}a
    + q^2d\beta\otimes d\cdot{\sf d}\gamma \right.\\
   &\quad -\left. q^2 ad\otimes\beta\cdot{\sf d}\gamma - \beta\gamma\otimes d\cdot{\sf d}a - q^2\beta d\otimes d\cdot{\sf d}\gamma
    - q^2a\gamma\otimes\beta\cdot{\sf d}a\right]\\
&= {\cal D}_q^{-1} \left[(\gamma\beta - q^2 ad)\otimes\beta\cdot{\sf d}\gamma + (da - \beta\gamma)\otimes d\cdot{\sf d}a \right] = {\bf 1}\otimes w_1.
\end{align*}
Similarly, one can see that the others are left-coinvariant. \hfill$\square$

From (\ref{4.23}), the defining formulae for $w_k$ can be inverted to
\begin{equation} \label{4.25}
{\sf d}t_{ij} = \sum_{k=1}^2 \,q^{p(t_{ik})} \, t_{ik} \cdot w_{kj}.
\end{equation}

\begin{remark} \label{rem4.7} 
The differentials of the generators of the algebra ${\cal O}(\widetilde{\rm GL}_q(2))$ satisfy the formula
\begin{align*}
{\sf d}f = q^{p(f)} \, \sum_{i=1}^4 \, (\chi_{_i}*f) \,w_i = q^{p(f)} \, \sum_{i=1}^4 \, m\circ({\rm id}\otimes\chi_{_i})\Delta(f) \,w_i.
\end{align*}
So, we have
\begin{align*}
\chi_{_a} &= \begin{pmatrix} 1 & 0 \\ 0 & 0 \end{pmatrix}, \quad \chi_{_\beta} = \begin{pmatrix} 0 & q \\ 0 & 0 \end{pmatrix}, \quad
\chi_{_\gamma} = \begin{pmatrix} 0 & 0 \\ q^2 & 0 \end{pmatrix}, \quad \chi_{_d} = \begin{pmatrix} 0 & 0 \\ 0 & 1 \end{pmatrix}.
\end{align*}
These matrices satisfy the following relations
\begin{align*}
[\chi_{_a},\chi_{_\beta}] = \chi_{_\beta}, \quad [\chi_{_a},\chi_{_\gamma}] = -\chi_{_\gamma}, \quad
[\chi_{_a},\chi_{_d}] = 0, \quad [\chi_{_\beta},\chi_{_\gamma}] = \chi_{_a} - \chi_{_d},
\end{align*}
where $[u,v]=uv-vu$.
\end{remark}

The following theorem determine the commutation relations between left-coinvariant forms and the elements of ${\cal O}(\widetilde{\rm GL}_q(2))$ for the Case I.

\begin{thm} \label{theo4.6} 
Let $(\Omega,\Delta_L)$ be a left-covariant bimodule over ${\cal O}(\widetilde{\rm GL}_q(2))$ and $\{w_k\}_{k=1}^4$ be a basis in the vector space of all left-invariant elements of $\Omega$. Then there exist linear map
$\mu:{\cal O}(\widetilde{\rm GL}_q(2))\longrightarrow M_4({\cal O}(\widetilde{\rm GL}_q(2)))$ such that
\begin{align} \label{4.26}
w_i \cdot f
&= q^{p(f)p(w_i)} \, \sum_{j=1}^4 \, \mu_{ij}(f)\, w_j
\end{align}
for all $f\in{\cal O}(\widetilde{\rm GL}_q(2))$.
\end{thm}

\noindent{\it Proof}.
Using (\ref{4.7}) and (\ref{4.15}), we find the matrices $\mu(f)$ for $f\in{\cal O}(\widetilde{\rm GL}_q(2))$ as
\begin{align} \label{4.27}
\mu(a) &= \left(\begin{matrix} P_1^{-1}a & 0 & (1-qP_1^{-1})\beta & 0\\ 0 & a & 0 & 0\\ 0 & 0 & qP_1a & 0\\
(qP_2^{-1}-1)a & 0 & (P_2^{-1}-q^2)\beta & a \end{matrix}\right), \nonumber\\
\mu(\beta) &=\left(\begin{matrix} qP_2^{-1}\beta & 0 & 0 & 0\\ 0 & P_2\beta & 0 & 0\\ 0 & 0 & q^2P_2^{-1}\beta & 0\\ 0 & (1-P_2^{-1})a & 0 & qP_2^{-1}\beta \end{matrix}\right), \nonumber\\
\mu(\gamma) &= \left(\begin{matrix} P_1^{-1}\gamma & 0 & (1-qP_1^{-1})d & 0\\ 0 & \gamma & 0& 0\\ 0 & 0 & qP_1\gamma & 0\\ (qP_2^{-1}-1)\gamma & 0 & (P_2^{-1}-q^2)d & \gamma \end{matrix}\right), \\
\mu(d) &= \left(\begin{matrix} qP_2^{-1}d & 0 & 0 & 0\\ 0 & P_2d & 0 & 0\\ 0 & 0 & q^2P_2^{-1}d & 0\\ 0 & (1-P_2^{-1})\gamma & 0 & qP_2^{-1}d \end{matrix}\right) \nonumber
\end{align}
where $P_1\in\{q,q^2\}$ and $P_2\in\{1,q^2\}$. Inserting these matrices in (\ref{4.26}) leads to commutation relations of Maurer-Cartan forms with the generators of the algebra
${\cal O}(\widetilde{\rm GL}_q(2))$. \hfill$\square$

\begin{thm} \label{theo4.7} 
The map $\mu$ is ${\mathbb C}$-linear ${\mathbb Z}_3$-graded homomorphism such that
\begin{equation} \label{4.28}
\mu_{ij}(fg) = \sum_{k=1}^4 q^{p(g)[p(w_k)-p(w_i)]} \,\mu_{ik}(f)\mu_{kj}(g), \quad \forall f,g\in{\cal O}(\widetilde{\rm GL}_q(2)).
\end{equation}
\end{thm}

\noindent{\it Proof}.
Using (\ref{4.26}), we can write
\begin{align*}
w_i\cdot(fg) &= q^{p(fg)p(w_i)} \, \sum_{j=1}^4 \mu_{ij}(fg) w_j.
\end{align*}
In the other hand, from the identity $w_i\cdot(fg) = (w_i\cdot f)g$,
\begin{align*}
q^{-p(f)p(w_i)} \, (w_i\cdot f)g
&= \sum_{j=1}^4 \mu_{ij}(f) w_j\cdot g = \sum_{j=1}^4 \mu_{ij}(f) \cdot q^{p(g)p(w_j)} \, \sum_{k=1}^4 \mu_{jk}(g) w_k \\
&= \sum_{j=1}^4\sum_{k=1}^4 q^{p(g)p(w_k)} \, \mu_{ik}(f) \mu_{kj}(g) w_j
\end{align*}
so that we conclude the equality in (\ref{4.28}). \hfill$\square$

\begin{remark} \label{rem4.8} 
Since we know, from \cite{Woro}, that
\begin{equation*}
\Delta\circ\mu = ({\rm id}\otimes\mu)\circ\Delta
\end{equation*}
we can write
\begin{equation*}
\mu(T) = T\cdot F(T)
\end{equation*}
in a compact form. Here $F(f)$ is a $4\times4$-matrix with entries in ${\mathbb C}$. Explicitly,
\begin{align} \label{4.29}
F(a) &= \left(\begin{matrix} P_1^{-1} & 0 & 0 & 0\\ 0 & 1 & 0 & 0\\ 0 & 0 & qP_1 & 0\\ qP_2^{-1}-1 & 0 & 0 & 1 \end{matrix}\right), \quad
F(\beta) = \left(\begin{matrix} 0 & 0 & 0 & 0\\ 0 & 0 & 0 & 0\\ 0 & 0 & 0 & 0\\ 0 & 1-P_2^{-1} & 0 & 0 \end{matrix}\right), \\
F(\gamma) &=\left(\begin{matrix} 0 & 0 & 1-qP_1^{-1} & 0\\ 0 & 0 & 0 & 0\\ 0 & 0 & 0 & 0\\ 0 & 0 & P_2^{-1}-q^2 & 0 \end{matrix}\right), \quad
F(d) = \left(\begin{matrix} qP_2^{-1} & 0 & 0 & 0\\ 0 & P_2 & 0 & 0\\ 0 & 0 & q^2P_2^{-1} & 0\\ 0 & 0 & 0 & qP_2^{-1} \end{matrix}\right) \nonumber
\end{align}
where $P_1\in\{q,q^2\}$ and $P_2\in\{1,q^2\}$. It is easy to see that the matrices in (\ref{4.29}) are a representation of $\widetilde{\rm GL}_q(2)$ with entries in ${\mathbb C}$, that is, they satisfy the relations (\ref{3.1}) in the sense that
\begin{equation*}
F_{ij}(f\cdot g) = \sum_{k=1}^4 F_{ik}(f)\cdot F_{kj}(g).
\end{equation*}
\end{remark}

We know that from \cite{Woro}, if $(\Omega,\Delta_L)$ is a left-covariant bimodule over an algebra ${\cal A}$ and $\{w_k\}_{k=1}^4$ is a basis in the vector space of all left-invariant elements of $\Omega$ then there exist linear functionals
$F_{ij}$ such that
\begin{equation*}
w_i \cdot f =
\sum_{j=1}^4 \,m\circ({\rm id}\otimes F_{ij})\Delta(f)\, w_j.
\end{equation*}
The ${\mathbb Z}_3$-graded form of this formula is expressed in the following corollary.

\begin{corollary} \label{cor4.5} 
Let $(\Omega,\Delta_L)$ be a left-covariant bimodule over ${\cal O}(\widetilde{\rm GL}_q(2))$ and $\{w_k\}_{k=1}^4$ be a basis in the vector space of all left-invariant elements of $\Omega$. Then there exist linear map
$F:{\cal O}(\widetilde{\rm GL}_q(2))\longrightarrow M_4({\mathbb C})$ such that
\begin{align} \label{4.30}
w_i \cdot f
&= q^{p(f)p(w_i)} \, \sum_{j=1}^4 \,m\circ({\rm id}\otimes F_{ij})\Delta(f)\, w_j
\end{align}
for all $f\in{\cal O}(\widetilde{\rm GL}_q(2))$.
\end{corollary}

\begin{remark} \label{rem4.9} 
There is a kind of inverse formula of (\ref{4.26}) as follows
\begin{align*}
f \cdot w_i
&= \sum_{j=1}^4 q^{-p(f)p(w_j)} \,  w_j \, \tilde{\mu}_{ji}(f), \quad \forall f\in{\cal O}(\widetilde{\rm GL}_q(2))
\end{align*}
where
\begin{equation*}
\tilde{\mu}_{ij}(fg) = \sum_{k=1}^4 q^{p(g)[p(w_i)-p(w_k)]} \,\tilde{\mu}_{ik}(f)\tilde{\mu}_{kj}(g)
\end{equation*}
for all $f,g\in{\cal O}(\widetilde{\rm GL}_q(2))$.
\end{remark}

\begin{remark} \label{rem4.10} 
The operators  $\mu$ and $\tilde{\mu}$ are related by the formulae
\begin{equation*}
\sum_{k=1}^4 q^{[p(f)-p(\mu_{ik}(f))]p(w_j)} \,\tilde{\mu}_{jk}(\mu_{ik}(f)) = f \, \delta_{ij}, \quad \forall f\in {\cal O}(\widetilde{\rm GL}_q(2)).
\end{equation*}
\end{remark}

Using (\ref{4.7}) and (\ref{4.15}) we find the commutation rules of the forms as follows.

\begin{thm} \label{theo4.8} 
The commutation rules of the forms are as follows
\begin{align} \label{4.31}
w_1 \cdot w_2 &= q^2P_1 \, w_2\cdot w_1, \quad w_1 \cdot w_3 = P_2 \, w_3\cdot w_1, \quad w_2 \cdot w_3 = P_1P_2^{-1} \, w_3\cdot w_2, \nonumber\\
w_1 \cdot w_4 &= q^2P_2 \, w_4\cdot w_1 + (1-q^2P_2) \, (w_1^2 - P_1P_2 \, w_3\cdot w_2), \nonumber\\
w_2 \cdot w_4 &= q \, w_4\cdot w_2 + q(qP_2^{-1}-1) \, w_1\cdot w_2), \\
w_3 \cdot w_4 &= q^2 \, w_4\cdot w_3 + (P_2^{-1}-q^2) \, w_1\cdot w_4, \nonumber
\end{align}
where $P_1\in\{q,q^2\}$ and $P_2\in\{1,q^2\}$.
\end{thm}

\subsection{The ${\mathbb Z}_3$-graded partial derivatives} \label{sec4.4}
In this section, we will continue with the relations given for the Case I in Theorem 4.5. To obtain the relations between the generators of ${\cal O}(\widetilde{\rm GL}_q(2))$ and derivatives, let us introduce the partial derivatives of the generators of the algebra.

Since $(\Omega,{\sf d})$ is a left covariant differential calculus, for any element $f$ in ${\cal O}(\widetilde{\rm GL}_q(2))$ there are uniquely determined elements $\partial_k(f)\in {\cal O}(\widetilde{\rm GL}_q(2))$ such that
\begin{equation}\label{4.32}
{\sf d}f = {\sf d}a \,\partial_a(f) + {\sf d}\beta \,\partial_\beta(f) + {\sf d}\gamma \,\partial_\gamma(f) + {\sf d}d \,\partial_d(f).
\end{equation}
For consistency, the grades of the derivatives $\partial_a$ and $\partial_d$ should be 0, and $\partial_\beta$ and $\partial_\gamma$ should be 1 and 2, respectively.

\begin{defin} \label{def4.6} 
The linear maps $\partial_a,\partial_d,\partial_\beta,\partial_\gamma:{\cal O}(\widetilde{\rm GL}_q(2))\longrightarrow{\cal O}(\widetilde{\rm GL}_q(2))$ defined by (\ref{4.32}) are called the {\it partial derivatives} of the calculus $(\Omega,{\sf d})$.
\end{defin}

The next theorem gives the relations between the elements of ${\cal O}(\widetilde{\rm GL}_q(2))$ and their partial derivatives.

\begin{thm} \label{theo4.9} 
The relations of the generators of the algebra ${\cal O}(\widetilde{\rm GL}_q(2))$ with partial derivatives are given by the formulae
\begin{equation}\label{4.33}
\partial_i(f\cdot g) = \partial_i(f)\cdot g + \sum_{j=1}^4 q^{-p(f)[1+p(\partial_j)]} \, \sigma_{ij}(f) \cdot \partial_j(g)
\end{equation}
for all $f,g \in {\cal O}(\widetilde{\rm GL}_q(2))$.
\end{thm}

\noindent{\it Proof}.
From the ${\mathbb Z}_3$-graded Leibniz formula for {\sf d},
\begin{equation*}
{\sf d}(f\cdot g) = ({\sf d}f)\cdot g + q^{p(f)} f\cdot ({\sf d}g),
\end{equation*}
and (\ref{4.9}) with (\ref{4.32}), we can obtain the twisted ${\mathbb Z}_3$-graded Leibniz formulas for partial derivatives. If
$f,g\in {\cal O}(\widetilde{\rm GL}_q(2))$, then
\begin{align*}
{\sf d}(f\cdot g)
&= \sum_{i=1}^4 {\sf d}f_i\cdot \partial_i(f\cdot g) = \left(\sum_{i=1}^4 {\sf d}f_i\cdot \partial_i(f)\right)\cdot g
    + q^{p(f)} f\cdot \sum_{i=1}^4 {\sf d}g_i\cdot \partial_i(g)\\
&=\sum_{i=1}^4 {\sf d}f_i\cdot \partial_i(f)\cdot g+q^{p(f)}\,\sum_{i=1}^4 q^{p(f)p({\sf d}g_i)}\,\sum_{j=1}^4{\sf d}g_j\cdot\sigma_{ji}(f)\cdot\partial_i(g) \\
&= \sum_{i=1}^4 {\sf d}f_i\cdot\left\{\partial_i(f)\cdot g + q^{p(f)[1+p({\sf d}f_i)]}\,\sum_{j=1}^4 \sigma_{ji}(f)\cdot\partial_i(g)\right\}
\end{align*}
so that we have the equality in (\ref{4.33}) by interchanging $i$ and $j$ in the last expression. \hfill$\square$

Since we have the condition(s) ${\sf d}^3=0$ (and ${\sf d}^2\ne0$) in ${\mathbb Z}_3$-graded calculus, it requires some tedious computations to find the relations of partial derivatives. Once the operator {\sf d} is applied to the equation (\ref{4.1}), one obtains
\begin{align*}
{\sf d}^2f
&= \sum_{i=1}^4 \left({\sf d}^2f_i \cdot\partial_if + q^{p({\sf d}f_i)} \,{\sf d}f_i \,\sum_{j=1}^4 {\sf d}f_j \cdot\partial_j(\partial_if)\right).
\end{align*}
Next, to obtain the desired relations the operator {\sf d} will be re-applied to both sides of this equation. As a result, we have the following theorem:

\begin{thm} \label{theo4.10}
The relations between partial derivatives are of the form
\begin{align} \label{4.34}
\partial_a \partial_\beta &= q^2P_2^{-1} \,\partial_\beta \partial_a, \quad \partial_a \partial_\gamma = \partial_\gamma \partial_a, \quad
   \partial_\beta \partial_d = q^2 \,\partial_d \partial_\beta, \quad \partial_\gamma \partial_d = q^2P_2^{-1} \,\partial_d \partial_\gamma,\nonumber\\
\partial_a \partial_d &= q^2P_2^{-1} \,\partial_d \partial_a, \quad \partial_\beta \partial_\gamma = q^2P_2^{-1}\, \partial_\gamma \partial_\beta + (1-q) \,\partial_d \partial_a,
\end{align}
where $P_2\in\{1,q^2\}$.
\end{thm}

\end{document}